\magnification=1200
\input amssym.def
\input amssym.tex
 \font\newrm =cmr10 at 24pt
\def\bul{\raise .9pt\hbox{\newrm .\kern-.105em } }

 \def\fr{\frak}

 \baselineskip 20pt
 
 \def\h{\hbox{ }}

 \def\hh{{\fr h}}

 \def\g{{\fr g}}

 \def\<{\le}
 \def\>{\ge}

 \def\s{{\h\subset\h}}
 
 \def\vs{\vskip }

 \def\mapright#1
  {\smash{\mathop
  {\longrightarrow}
  \limits^{#1}}}

 \def\kk#1{{\kern .4 em} #1}
 \def\vs{\vskip 1pc}

\font\authorfont=cmr10 at 12pt
\font\titlefont=cmr12 at 14pt

\font\ninerm=cmr10 at 9pt
	
\rm

\hsize = 30pc
\vsize = 45pc
\overfullrule = 0pt

\rm
\centerline{\titlefont The Coxeter element and the branching law for the }
\centerline{\titlefont finite subgroups of $SU(2)$}\vskip 2pc 
 \baselineskip=11pt
\vskip8pt
\centerline{\authorfont BERTRAM
KOSTANT\footnote*{\ninerm Research supported in part by NSF grant
DMS-0209473 and in part by the \hfil\break KG\&G Foundation.}}\vskip 2pc	
\rm
\baselineskip 20pt

	\centerline{\bf 0. Introduction}\vskip 1.5pc 0.1. Let $\Gamma$ be a finite subgroup of $SU(2)$. The question we
will deal with in this paper is how an arbitrary (unitary) irreducible representation of $SU(2)$ decomposes
under the action of $\Gamma$. The theory of McKay assigns to $\Gamma$ a complex simple Lie algebra $\g$ of type
$A-D-E$. The assignment is such that if $\widetilde {\Gamma}$ is the unitary dual of $\Gamma$ we may parameterize
$\widetilde {\Gamma}$ by the nodes (or vertices) of the extended Coxeter-Dynkin diagram of
$\g$.

Let $\ell = rank\,\g$ and let $I = \{1,\ldots,\ell\}$. Let $I_{ext} = I\cup \{0\}$. The nodes may be
identified with a set of simple roots of the affine Kac-Moody Lie algebra associated to $\g$ and are indexed by
$I_{ext}$. We can then write $\Gamma = \{\gamma_i\},\,i\in I_{ext}$. Let $\Pi = \{\alpha_i\},\,i\in I$, be
the set of simple roots of $\g$ itself. One has $\gamma_0$ is the trivial 1 dimensional representation of
$\Gamma$ and, for $i\in I$, $$dim\,\gamma_i = d_i\eqno (0.1)$$ where $$\psi = \sum_{i\in I} d_i\,\alpha_i\eqno
(0.2)$$ is the highest root. For proofs and details about the McKay correspondence see e.g. [G-S,V], [M] and
[St]. \vs 0.2. The unitary dual of $SU(2)$ is indexed by the set $\Bbb Z_+$ of nonnegative integers and will be
written as
$\{\pi_n\},\,n\in \Bbb Z_+$ where $$dim\,\pi_n = n+1\eqno (0.3)$$ Our problem is the determination of $m_{n,i}$
where $n\in \Bbb Z_+,\,\,i\in I_{ext}$ and $$m_{n,i} = \hbox{multiplicity of $\gamma_i$ in $\pi_n|\Gamma$}$$ It
is resolved with the determination of the formal power series $$m(t)_i = \sum_{n=0}^{\infty}m_{n,i}\,t^n\eqno
(0.4)$$ To do this one readily notes that it suffices to consider only the case where $\Gamma = F^*$ and $F^*$
is the pullback to
$SU(2)$ of a finite subgroup $F$ of $SO(3)$. This eliminates only the case where $\Gamma$ is a cyclic group of odd
order and $\g$ is of type $A_{\ell}$ where $\ell$ is even. For the remaining cases the Coxeter number $h$ of $\g$
is even and we will put $$g = h/2\eqno (0.5)$$ Also for the remaining cases there is a special index $i_*\in I$.
If $\g$ is of type $D$ or $E$ then $\alpha_{i_*}$ is the branch point of the Coxeter-Dynkin diagram of $\g$. 
If $\g$ is of type $A_{\ell}$ then $\alpha_{i_*}$ is the midpoint of the diagram (recalling that $\ell$ is odd). 

If $i=0$ the determination of $m(t)_0$ is
classical and is known from the theory of Kleinian
singularities. In fact there exists positive integers
$a<b$ such that $$m(t)_0 = {1 + t^h\over (1-t^a)(1-t^b)}\eqno
(0.6)$$ The numbers $a$ and $b$ in Lie theoretic terms is given in 
 \vs {\bf Theorem 0.1}. {\it One has $a = 2d_{i_*}$ and $b$ is given by the condition
that $$ \eqalign{a\,b &= 2\,|F^*|\cr &= 4\,|F|\cr}$$}\vs It remains then to determine $m(t)_i$ for $i\in I$.
\vs {\bf Proposition 0.2.} {\it If $i\in I$ there exists a polynomial $z(t)_i$ of degree less than $h$ such that
$$m(t)_i = {z(t)_i\over (1-t^a)(1-t^b)}\eqno (0.7)$$}\vs The problem then is to determine the polynomial
$z(t)_i$. This problem was solved in [K] using the orbits of a Coxeter element $\sigma$ on a set of roots $\Delta$
for
$\g$. In the present paper we will put the main result of [K] in a simplified form. See Theorem 1.13 in the
present paper. Also the present paper makes explicit some results that are only implicit in [K]. For example
introducing $\widetilde {\Pi}$  (see (1.10)) and making the assertions in Remark 10 and Theorems 8, 9, 11 and 12. 

The set $\Pi$ generates a system, $\Delta_+$, of posiive roots. The highest root
$\psi\in
\Delta$ defines a certain subset $\Phi\s \Delta_+$ of cardinality $2h-3$. Because of its connection with a
Heisenberg subalgebra of $\g$ this subset is referred to as the Heisenberg subsystem of $\Delta_+$. The new
formulation explicitly shows how the polynomials $z(t)_i$ arise from the intersection $$(\hbox{orbits of the
Coxeter element $\sigma$})\cap (\hbox{the Heisenberg subsystem}\,\,\Phi)\eqno (0.8)$$

The polynomials $z(t)_i$ are listed in [K]. The special case where $\g$ is of type $E_8$ is also given
in the present paper (see Example 1.7.). Unrelated to the Coxeter element the polynomials
$z(t)_i$ are also determined in Springer, [Sp]. They also appear in another context in Lusztig, [L1] and [L2].
Recently, in a beautiful result, Rossmann, [R], relates the character of $\gamma_i$ to the polynomial
$z(t)_i$.\vskip 1.5pc \centerline{\bf 1. The main result - Theorem 1.13.}\vskip 1.5pc 1.1. Proofs of the main
results stated here are in [K]. 

Let $F$ be a finite subgroup of $SO(3)$ and let $$F^*\s SU(2)\eqno
(1.1)$$ be the pullback of the double covering $$SU(2)\to SO(3)$$ 
The unitary dual $\widehat {SU(2)}$ of $SU(2)$ is represented by the
set
$\{\pi_n\},\,n\in \Bbb Z_+$, where if $S(\Bbb C^2)$ is the symmetric
algebra then $$\pi_n:SU(2)\to S^n(\Bbb C^2)$$ is the
$n+1$ dimensional representation defined by the natural action of
$SU(2)$ on $\Bbb C^2$. We are ultimately interested 
in determining
 how the restriction $\pi_n|F^*$ decomposes into irreducible
representations of the finite subgroup $F^*$, for any $n$, and
relating this determination to the structure of the simple Lie
algebra corresponding to $F^*$ by the McKay correspondence. We now
recall this correspondence. 

Let $\g$ be a complex simple Lie algebra and let $\hh$ be a
Cartan subalgebra of $\g$. Let $\ell = rank\,\g$, and if $\hh'$ is
the dual space to $\hh$, let $\Delta\s \hh'$ be the set of roots
for the pair $(\g,\hh)$. Let $W$, operating in $\hh'$, be the Weyl group. Let
$\Pi$ be the set of simple positive roots with respect to a choice, $\Delta_+$, of 
positive roots. If $I = \{1,\dots, \ell\}$ we will write $\Pi =
\{\alpha_i\},\,i\in I$. We may regard $\Pi$ as the vertices (or nodes) of the
Coxeter-Dynkin diagram associated to $\g$.  The extended
Coxeter-Dynkin diagram has an additional node $\alpha_0$. 

The McKay
correspondence assigns to
$F^*$ a complex simple Lie algebra $\g = \mu(F^*)$ of type $A-D-E$. 
The assignment has a number of properties: (1), the
unitary dual
$\widehat{F^*}$ may be parameterized by indices of the nodes of
the extended Coxeter-Dynkin diagram of $\g$. In particular
$card\,\widehat{F^*} = \ell + 1$ and we can write
$\widehat{F^*}=\{\gamma_i\},\,i\in I_{ext} = I\cup\{0\}$. Next
(2), $\gamma_0$ is the trivial 1-dimensional representation, and if
$i\in I$, then $$dim\,\gamma_i = d_i$$ where $$\psi =
\sum_{i=1}^{\ell} d_i\alpha_i$$ is the highest root in $\Delta$.
In addition (3), if $\gamma$ is the two-dimensional representation
defined by (1.1) and $A$ is the $(\ell+1)\times (\ell+1)$ matrix
defined so that $$\gamma_i\otimes \gamma = \sum_{j=0}^{\ell}
A_{ij}\gamma_j\eqno (1.2)$$ then $C$ is the Cartan matrix of the
extended Coxeter-Dynkin of $\g$ where $$ C_{ij} = 2\delta_{ij}-
A_{ij}$$ \vskip .5pc
1.2. Returning to our main problem, for $i\in I_{ext}$ and $n\in \Bbb
Z_+$, let
$$m_{n,i} =\hbox{multiplicity of $\gamma_i$ in $\pi_n|F^*$}$$ and
introduce the generating formal power series $$m(t)_i =
\sum_{n=0}^{\infty}\,m_{n,i}\,t^n$$ If $i=0$, the determination of $m(t)_i$ is
classical and is known from the theory of Kleinian
singularities. That is, in this case $$m_{n,0}= dim\,(S^n(\Bbb
C))^{F^*}$$ In fact let $h$ be the Coxeter number of $\g$ so that
$$\ell( h+1) = dim\,\g$$ Then there exists positive integers
$a<b$ such that $$m(t)_0 = {1 + t^h\over (1-t^a)(1-t^b)}\eqno
(1.3)$$ To define the numbers $a$ and $b$ in Lie theoretic terms
one notes that $\mu(F^*)$ is of type $D$, $E$ or $A_{\ell}$ where
$\ell$ is odd. In any of these case there is a special index $i_*\in
I$. If $\mu(F^*)$ is of type $D$ or $E$, then $\alpha_{i_*}$ is the
branch point of the Coxeter-Dynkin diagram of $\g$. If $\mu(F^*)$
is of type $A_{\ell}$, then $\alpha_{i_*}$ is the midpoint of the
diagram (recall that $\ell$ is odd in this case). \vs {\bf Theorem
1.1}. {\it One has $a = 2d_{i_*}$ and $b$ is given by the condition
that $$ \eqalign{a\,b &= 2\,|F^*|\cr &= 4\,|F|\cr}\eqno
(1.4)$$}\vskip .5pc See Lemma 5.14 in [K]. The cases under consideration are
characterized by the condition that $h$ is even. We put $g =
h/2$. The parity of $g$ will play a later
role. \vs {\bf Remark 1.2.} One proves (see Lemma 5.7 in [K]) that $b$ may also be given by
$$ b = h+2-a\eqno (1.5)$$ so that $b$, as well as $a$, is even.\vs
 The following table lists the various cases under consideration.
In the table
$\Delta_n$ is the dihedral group of order $2n$. 

$$\matrix{\underline {F}&\underline{\g}&& \underline
{a}& &\underline {b}& &\underline {h}&& \underline{g}\cr & & & & & &
& & &\cr
\Bbb Z_n&A_{2n-1}&&2 && 2n && 2n&& n\cr \Delta_n& D_{n+2}&&4&&
2n&&2n+2&&n+1\cr Alt_4&E_6&&6&&8&&12&&6\cr
Sym_4&E_7&&8&&12&&18&&9\cr Alt_5&E_8&& 12&&20&& 30&&15\cr}$$\vskip
1pc {\bf Proposition 1.3.} {\it There exists a unique partition $$\Pi
= \Pi_1\cup \Pi_2\eqno (1.6)$$ such that if $k=1,2$ and
$\alpha_i,\alpha_j,\in \Pi_k$ where $i\neq j$ then $\alpha_i$ is
orthogonal to $\alpha_j$. Furthermore all the roots in $\Pi_2$ are
orthogonal to the highest root $\psi$, or equivalently the root
$\alpha_0$ is orthogonal to all the roots in $\Pi_2$.}\vs One has the
disjoint union $I = I_1\cup I_2$ where, if $k\in \{1,2\}$, $\Pi_k =
\{\alpha_i\mid i\in I_k\}$.\vs  {\bf Remark 1.4.} It is immediate from
(1.2) that if
$A_{ij}\neq 0$ and
$i\in I_k$ then $j$ is in the complement of $I_k$ in $I$. It
then follows that $\gamma_i$ descends to a representation of $F$
(i.e., $\gamma_i(-1) = 1$) if and only if $k= 2$. In particular
$$m_{n,i} = 0\,\,\,\hbox{if $n$ and $k$ have opposite parities
where $\alpha_i\in \Pi_k$.}\eqno (1.7)$$  \vs If $i\in I$ let
$s_i\in W$ be the reflection defined by $\alpha_i$ so that $s_i$
commutes with $s_j$ if $i,j\in I_k,\,k\in \{1,2\}$. Put $\tau_k
=\prod_{i\in I_k}\,s_i$. Then $$\eqalign{\tau_1^2 &= \tau_2^2\cr
&=identity\cr}$$ One defines a Coxeter element $\sigma\in W$  by
putting
$$\sigma = \tau_2\,\tau_1\eqno (1.8)$$\vskip .5pc {\bf Remark 1.5.} 
Every element in $W$ is contained in a dihedral subgroup of $W$.
Since, as one knows, the centralizer of a Coxeter element is the
cyclic group (necessarily of order $h$) generated by the Coxeter
element, a dihedral group containing the Coxeter element is unique.
It is clear that $\tau_1$ and $\tau_2$ are in the dihedral group
containing $\sigma$ and, in fact, are in the complementary coset of the
cyclic group generated by $\sigma$.\vs As a extension of (1.3)
one knows (see (5.7.2) in [K]) that for any $i\in I$ there exists a polynomial
$z(t)_i$ of degree less than $h$ such that $$m(t)_i = {z(t)_i\over
(1-t^a)(1-t^b)}\eqno (1.9)$$ so that $m(t)_i$ is known as soon as
one knows the polynomial $z(t)_i$.\vs {\bf Remark 1.6.} Note that
by (1.6) and evenness of $a$ and $b$ (Remark 1.2) one must have
that the only powers of $t$ which have a nonzero coefficient are odd if
$i\in I_1$ and even if $i\in I_2$.\vskip 1.5pc
{\bf Example 1.7}. Consider the case
where
$F$ is the icosahedral group so that
$\mu(F^*) = E_8$. In the listing of $z(t)_i$ below we will replace the
arbitrary index $i$ by the more informative 
$\{d_i\}$. Since there exists in certain cases two distinct
$i,j\in I$ such that $dim\,\gamma_i = dim\,\gamma_j$ we will write
$\underline {\{d_j\}}$ for $j$ when the ``distance" of $\alpha_j$ to
$\alpha_0$ is greater than the ``distance" of $\alpha_i$ to
$\alpha_0$. Note that $d_{i_*} = 6$. \vskip 1.5pc
\item \qquad\qquad ${z(t)_{\{2\}}} = t + t^{11} + t^{19} + t^{29}$
\item \qquad\qquad ${z(t)_{\{3\}}} = t^2 + t^{10} + t^{12} + t^{18} +
t^{20} + t^{28}$
\item \qquad\qquad ${z(t)_{\{4\}}} = t^3 + t^{9} + t^{11} + t^{13} +
t^{17} + t^{19} + t^{21} + t^{27}$
\item \qquad\qquad ${z(t)_{\{5\}}} = t^4 + t^{8} + t^{10} + t^{12} +
 t^{14} + t^{16} + t^{18} + t^{20} + t^{22} + t^{26}$
\item \qquad\qquad ${z(t)_{\{6\}}} = t^5 + t^{7} + t^{9} + t^{11} + t^{13}
+ 2\,t^{15} + t^{17} + t^{19} + t^{21} + t^{23} + t^{25}$
\item \qquad\qquad ${z(t)_{\underline{\{4\}}}} = t^6 + t^8 + t^{12} +
t^{14} + t^{16} + t^{18} + t^{22} + t^{24}$
\item \qquad\qquad ${z(t)_{\underline{\{2\}}}} = t^7 + t^{13} + t^{17} +
t^{23}$
\item \qquad\qquad ${z(t)_{\underline{\{3\}}}} = t^6 + t^{10} + t^{14} +
t^{16} + t^{20} + t^{24}$\vskip 1.5pc  We now modify $\Pi$ by defining
$$\widetilde {\Pi} =\{\beta_i\mid i\in I\}\eqno (1.10)$$  where
$\beta_i = \alpha_i$ if $i\in I_1$ and $\beta_i = -\alpha_i$ if $i\in
I_2$. Let $Z\s W$ be the cyclic group generated by the Coxeter
element $\sigma$. Recall $(h+1)\,\ell$ so that $$card\,\Delta=
h\,\ell\eqno (1.11)$$ We have shown that $\sigma$ has $\ell$ orbits in
$\Delta$, each with $h$-elements, and that each orbit contains a unique element of $\widetilde {\Pi}$. That is,
one has 
\vs {\bf Theorem 1.8.} {\it For any $i\in I$ the
$\sigma$-orbit $Z\cdot\beta_i$ has $h$ elements and one has the
disjoint union $$\Delta =\sqcup_{i=1}^{\ell}\,\, Z\cdot \beta_i\eqno
(1.12)$$}\vs This result is readily proved using (6.9.2) in [K]. 

For any $i\in I$ let $(Z\cdot \beta_i)_+ = \Delta_+\cap Z\cdot
\beta_i$. One has (see (0.5)) $$\Delta_+ = g\,\ell\eqno (1.13)$$ \vskip .5pc {\bf
Theorem 1.9.} {\it For any $i\in I$ one has $card\,(Z\cdot \beta_i)_+ =
g$ and the disjoint union $$\Delta_+ = \sqcup_{i\in I}\,(Z\cdot
\beta_i)_+\eqno (1.14)$$}\vs It follows from (5.6.2) in [K] that (see (0.5)) 
$$\alpha_{i_*}\in \Pi_2\,\,\hbox{if $g$ is even and}\,\, \alpha_{i_*}\in
\Pi_1\,\,\hbox{if $g$ is odd.}\eqno (1.15)$$ \vskip .5pc Let
$\kappa$ be the long element of the Weyl group. One has (see Lemma 4.9 in
[K]) the following result of Steinberg: $$\sigma^g = \kappa\eqno (1.16)$$
so that
$\kappa\in Z$. \vs {\bf Remark 1.10.} Recall that
$\psi$ is the highest root. It is a consequence of
(5.6.2) in [K] that one has $\psi$ and $\beta_{i_*}$ are in the same $\sigma$
orbit. In fact if $g$ is odd then $$\eqalign{\sigma^{g-1\over 2}(\psi)
&=\beta_{i_*}\cr &= \alpha_{i_*}\cr}\eqno (1.17) $$ and if $g$ is even
then $$\eqalign{\sigma^{g\over 2}(\psi)
&=\beta_{i_*}\cr &= -\alpha_{i_*}\cr}\eqno (1.18)$$ One easily has that
$\sigma^g$ commutes with $\tau_1$ and $\tau_2$ so that, for $k\in
\{1,2\},$
$$\sigma^g(\Pi_k) = -(\Pi_k)\eqno (1.19)$$ Furthermore since $\kappa
(\psi) = -\psi$ one has that $$\sigma^g(\alpha_{i_*}) =
-\alpha_{i_*}\eqno (1.20)$$ so that in any case
$$\psi\,\,\hbox{and}\,\,\alpha_{i_*}\,\,\hbox{lie in the same
$\sigma$-orbit}\eqno (1.21)$$\vs 1.3. We come now to the main result---the
determination of $z(t)_i$ in terms of the orbit structure of $\sigma$
on $\Delta$. For any $\varphi\in \Delta_+$ let $i_{\varphi}\in I$ be
defined so that (by Theorem 1.9) $$\varphi\in (Z\cdot
\beta_{i_{\varphi}})_+\eqno (1.22)$$ But then there exists
$k_{\varphi}\in \{1,2\}$ such that $$i_{\varphi}\in
I_{k_{\varphi}}\eqno (1.23)$$ The following result follows from
(6.9.2) in [K]. \vskip .5pc
{\bf Theorem 1.11.} {\it Let
$\varphi\in \Delta_+$. Then there exists a unique positive integer
$n(\varphi)$ where $1\leq n(\varphi)\leq h$ with the same parity as
$k_{\varphi}$ such that if
$k_{\varphi}=1$ then $$\sigma^{n(\varphi)-1\over 2}(\varphi) =
\beta_{i_{\varphi}}\eqno (1.24)$$ If $k_{\varphi}= 2$
then $$\sigma^{n(\varphi)\over 2}(\varphi) =
\beta_{i_{\varphi}}\eqno (1.25)$$ }\vs One also has (see Remark 6.10 in
[K]) \vs {\bf Theorem 1.12.} {\it For any $i\in I_1$ the map
$$(Z\cdot \beta_i)_+\to \{0,1,\ldots,g-1\},\qquad \varphi \mapsto
{n(\varphi)-1\over 2}\eqno (1.26)$$ is a bijection and for any $i\in
I_2$ the map 
$$(Z\cdot \beta_i)_+\to \{1,\ldots,g\},\qquad \varphi \mapsto
{n(\varphi)\over 2}\eqno (1.27)$$ is a bijection.}\vs Let
$(\varphi,\varphi')$ be the restriction to $\Delta$ of the
$W$-invariant bilinear form on $\hh'$ induced by the Killing form on
$\g$. Let $\Phi =
\{\varphi\in \Delta\mid (\psi,\varphi) >0\}$. One easily has that
$\Phi\s\Delta_+$. Obviously $\psi\in \Phi$. One has $$card\,\Phi =
2\,h-3\eqno (1.28)$$ Because of its connection with a Heisenberg
subalgebra of $\g$ we refer to $\Phi$ as the Heisenberg subsystem of
$\Delta_+$. For $i\in I$ let $\Phi^i = \Phi\cap (Z\cdot \beta_i)_+$.
Our main result is \vs {\bf Theorem 1.13.} {\it Let $i\in I-\{i_*\}$. Then
$$z(t)_i = \sum_{\varphi\in \Phi^i}\,t^{n(\varphi)}\eqno (1.29)$$
Furthermore $$card\,\Phi^i = 2d_i\eqno (1.30)$$ In addition all the 
coefficients of
$z(t)_i$ are either $1$ or $0$ so that $$z(1)_i = 2\,d_i \eqno (1.31)$$
For $i=i_*$ one has $$z(t)_{i_*} = 2\,t^g + \sum_{\varphi\in
\Phi^{i_*},\,\varphi \neq \psi} t^{n(\varphi)}\eqno (1.32)$$ In
addition the coefficient of $t^g$ is 2 and all the other coefficients of
$z(t)_{i_*}$ are either $0$ or $1$. One also has $$\eqalign{z(1)_{i_*} &=
2\,d_{i_*}\cr &=a\cr}\eqno (1.33)$$ Finally $$z(t)_{i_*} = t^{g -a +2} +
t^{g-a +4} +
\cdots + t^{g-2} + 2\,t^g + t^{g +2} +\cdots + t^{g + a -4} + t^{g
+a-2}\eqno (1.34)$$} \vs Theorem 1.13 combines Theorem 6.6 and Lemma 6.14 in [K]. We note also that the 
expression (1.32) for $z(t)_{i_*}$ in Theorem 1.13 follows from the proof of Theorem 6.6 in [K]
 (see especially (5.8.1) in [K]). 
\vskip 1pc

\centerline{\bf References}\vskip 1.5pc
\rm
\item {[G-S,V]} Construction g\'eom\'etrique de la correspondance
de McKay. {\it Ann. Sci. Ecole Norm. Sup.} {\bf 16}, n$^o$3
(1983), 410-449.
\item {[K]} B. Kostant, The McKay correspondence, the Coxeter
element and representation theory. In {\it \'Elie Cartan
et les math\'ematiques d'aujourd'hui} (Lyon 1984). Ast\'erisque,
hors s\'erie, (1985), 209-255.
\item {[L1]}  G. Lusztig, Some examples of square integrable
representations of semisimple p-adic groups, {\it Trans. AMS} {\bf
277}(1983), 153-215.
\item {[L2]}  G. Lusztig, Subregular nilpotent elements and bases in
K-theory. {\it Canad. J. Math.}, {\bf 51} (6) (1999),
1194-1225.
\item {[M]} J. McKay, Graphs, singularities and finite
groups. {\it Proc. Symp. Pure math.}, {\bf 37} (1980), 183-186.
 \item {[R]} W. Rossmann, McKay's correspondence and characters of
finite subgroups of SU(2), {\it Noncommutative Harmonic Analysis,
in honor of Jacques Carmona}, Prog.in Math. {\bf 220}(2004),
Birkh\"auser, 441-458.
\item {[Sp]} T. Springer, Poincar\'e series of binary polyhedral
groups and McKay's correspondence. {\it Math. Ann.} {\bf
278}(1985), 587-598. 
\item {[St]} Finite subgroups of SU$_2$, affine Dynkin diagrams and
affine Coxeter elememts. Pac. J. Math. {\bf 118}(1985), 587-598,
Preprint 1982.

\parindent=30pt
\baselineskip=14pt
\vskip 1.9pc
\vbox to 60pt{\hbox{Bertram Kostant}
      \hbox{Dept. of Math.}
      \hbox{MIT}
      \hbox{Cambridge, MA 02139}} \noindent E-mail
kostant@math.mit.edu\vskip 1pc
 
\end

\centerline{\bf References}\vskip 1.5pc
\rm
\item {[G-S,V]} Construction g\'eom\'etrique de la correspondence
de McKay. {\it Ann. Sci. Ecole Nor. Sup.} {\bf 16}, n$^o$3
(1983), 410-449.
\item {[K]} B. Kostant, The McKay correspondence, the Coxeter
element and representation theory. In {\it \'Elie Cartan
et les math\'ematique d'aujourd'hui} (Lyon 1984). Ast\'erisque,
hors s\'erie, (1985), 209-255.
\item {[L1]}  G. Lusztig, Some examples of square integrable
representations of semisimple p-adic groups, {\it Trans. AMS} {\bf
277}(1983), 153-215.
\item {[L2]}  G. Lusztig, Subregular nilpotent elements and bases in
K-theory. {\it Canad. J. Math.},{\bf 51} (6) (1999),
1194-1225.
\item {[M]} J. McKay, Graphs, singularities and finite
groups. {\it Proc. Symp. Pure math.}, {\bf 37}(1980), 183-186.
 \item {[R]} W. Rossmann, McKay's correspondence and characters of
finite subgroups of SU(2), {\it Noncommutative Harmonic Analysis,
in honor of Jacques Carmona}, Prog. in Math. {\bf 220}(2004),
Birkhauser, 441-458.
\item {[Sp]} T. Springer, Poincar\'e series of binary polyhedral
groups and McKay's correspondence. {\it Math. Ann.} {\bf
278}(1985), 587-598. 
\item {[St]} Finite subgroups of SU$_2$, affine Dynkin diagrams and
affine Coxeter elememts. Pac. J. Math. {\bf 118}(1985), 587-598,
Preprint 1982
\parindent=30pt
\baselineskip=14pt
\vskip 1.9pc
\vbox to 60pt{\hbox{Bertram Kostant}
      \hbox{Dept. of Math.}
      \hbox{MIT}
      \hbox{Cambridge, MA 02139}} \noindent E-mail
kostant@math.mit.edu\vskip 1pc